# Estimating the Causal Effects of Marketing Interventions Using Propensity Score Methodology

Donald B. Rubin and Richard P. Waterman

*Abstract.* Propensity score methods were proposed by Rosenbaum and Rubin [*Biometrika* **70** (1983) 41–55] as central tools to help assess the causal effects of interventions. Since their introduction more than two decades ago, they have found wide application in a variety of areas, including medical research, economics, epidemiology and education, especially in those situations where randomized experiments are either difficult to perform, or raise ethical questions, or would require extensive delays before answers could be obtained. In the past few years, the number of published applications using propensity score methods to evaluate medical and epidemiological interventions has increased dramatically. Nevertheless, thus far, we believe that there have been few applications of propensity score methods to evaluate marketing interventions (e.g., advertising, promotions), where the tradition is to use generally inappropriate techniques, which focus on the prediction of an outcome from background characteristics and an indicator for the intervention using statistical tools such as least-squares regression, data mining, and so on. With these techniques, an estimated parameter in the model is used to estimate some global "causal" effect. This practice can generate grossly incorrect answers that can be self-perpetuating: polishing the Ferraris rather than the Jeeps "causes" them to continue to win more races than the Jeeps ⇔ visiting the high-prescribing doctors rather than the low-prescribing doctors "causes" them to continue to write more prescriptions. This presentation will take "causality" seriously, not just as a casual concept implying some predictive association in a data set, and will illustrate why propensity score methods are generally superior in practice to the standard predictive approaches for estimating causal effects.

*Donald B. Rubin is John L. Loeb Professor of Statistics, Harvard University, Cambridge, Massachusetts 02138, USA e-mail: rubin@stat.harvard.edu. Richard P. Waterman is Adjunct Associate Professor, The Wharton School, University of Pennsylvania, Philadelphia, Pennsylvania 19104, USA e-mail: waterman@wharton.upenn.edu.*

## 1. INTRODUCTION

This presentation is very simple in some sense, but in our experience the issues being discussed are often misunderstood, despite their importance. The application that is used throughout the first sections is a real one that was encountered nearly a decade ago, but we think has enough in common with many applications in business to be of general interest. The basic situation involves ordering a list of individuals to contact (e.g., by telephone or personal visit), from most likely to generate additional revenue to least likely to do so. Every effort at contact requires some investment, and we want to target for contact those individuals who are most likely to generate a return on that investment (ROI). The basic confusion in this situation is due to the more general confusion between a "before–after" change associated with an intervening event and a change that is "caused by" that intervening event; the culprit here is the word "change"—change from what? The before–after comparison is a change in time from before to after, but is that a "causal change"? The answer is almost always "no."

Our specific application was a project for a major pharmaceutical company concerned with their marketing interventions with doctors for the purpose of promoting sales of a particular "life-style" drug. The marketing interventions could be visiting a doctor to describe the details of the drug (so called, "detailing"), or it could be dining the doctor at a nice restaurant to convey similar information, or it could be providing free samples of the drug. All of these interventions, and other similar ones, are designed to lead to "more" prescriptions (scripts) for the drug written by the detailed doctor. But the critical question is: "more" than what? The answer is quite clear: more than that doctor would have written without the visit, dinner or free sample. Otherwise, the investment has had no positive return. Marketing interventions are designed to CAUSE A DIFFERENCE, and this difference, or change, is generally NOT a change in time.

The causal effect of the intervention on a doctor is the comparison of something you can see (e.g., the number of scripts written after being visited) with something you cannot see (e.g., the number of scripts written during the same period of time without being visited). Causal effects can be well estimated by essentially no existing statistical software based on predictive approaches, because causal effect estimation differs from simple prediction. Nevertheless, causal effects can often be well estimated in examples like ours, by propensity score technology, described and illustrated later, starting in Section 3. Essentially, the idea is to create matched pairs of units, where one member of the pair has been exposed to the intervention and the other has not, but they are otherwise identical before the time of exposure, that is, they are "clones." Finding such clones is a tall order because exact matches are almost impossible to find in realistically sized data sets, and this is where the propensity score technology enters.

In the next section we describe the difference between (1) simple prediction from the past to the future and (2) causal effects, and illustrate this distinction in a couple of totally trivial, but hopefully revealing, artificial examples. We then describe in Section 3 the idea of "cloning" for causal effect estimation—not a new idea, but hopefully expressed in such a way that makes important points transparent. Here we also introduce propensity score techniques.

The real example that motivated this presentation will then be described in Section 4, and diagnostic information will be presented concerning how successful the cloning using propensity scores appeared to be in this example. The results of our approach are estimates of individual doctor-level causal effects, which could be used as building blocks for addressing complex causal questions involving ROI. Specifically, these estimated doctor-level causal effects were then used to create an ordered list of the unvisited doctors, summarized in Section 5, ranked from those having large estimated causal effects of a visit, who should be visited, to those having small estimated causal effects of a visit, who should not be visited. In Section 6 we present an evaluation of our ordered list versus the company's standard (or traditional) ordering, and document the superiority of our causal ordering over their standard ordering, using the company's own criteria based on future scripts written.

Section 7 presents a more mechanical and general description of the basic methodology (e.g., in terms of units of analysis rather than doctors). Section 8 continues with a brief description of possible opportunities for the general approach in e-commerce. Section 9 concludes with a discussion of three key features of our general approach: the absence of any outcome variables when creating the clones; the opportunity to refine the causal estimates using models relating the outcome variables to background characteristics; and the use of traditional prediction models to select units, as defined by their observed background variables, that can be anticipated to have large causal effects of the intervention, and thus, a large ROI.

## 2. A CAUSAL EFFECT IS A "CHANGE," BUT NOT A CHANGE IN TIME

Display 1 is a very simple display of the title of this section. We have one doctor, and at time 1, we have in the left box the number of scripts that doctor has written in the six months prior to time 1. We have to make a choice to visit this doctor to provide details about the drug of interest, or not to visit. The top branch of the display represents what will happen if we visit, that is, detail, the doctor, where the box at the upper right gives the number of scripts written in the six months following the visit, up until time 2. In contrast, the bottom branch represents what will happen if we do not visit this doctor, and the box at the bottom right gives the number of scripts written during the same period of time if the doctor is not detailed.

The number of scripts written at time 2 given in the upper right box compared to the number at time 1 given in the left box is a change in time of the number of scripts written, but it is not the causal effect of the visit on number of scripts. It is the change in scripts written from time 1 to time 2 when the doctor is visited in between. Analogously, the number of scripts in the lower right box compared to the number of scripts in the left box also is a change of scripts written, and it is also a change in time, but is not the causal effect of not being visited on the number of scripts written.

The critical comparison here that is causal is the comparison of the number of scripts in the top right box and the number of scripts in the bottom right box, which is the causal effect of the doctor being visited versus not visited on the number of scripts written, which does not involve the box on the left at all, at least not without some overly strong assumption (e.g., the time-2 box without the visit is identical to the time-1 box—no change in time if not visited).

A causal effect is the comparison of the outcome that would be observed with the intervention and without the intervention, both measured at the same point in time. This is indicated by the comparison of apples with apples at time 2, whereas any comparison of something at time 2 with something at time 1 is indicated by the comparison of apples with oranges. This point, we know, is obvious, but its force is sometimes lost in the complication of real and hypothetical examples. The basic framework is often described as the "Rubin Causal Model" (Holland, 1986) for a sequence of articles starting in the 1970s, although the ideas obviously have much older roots (e.g., see Rubin, 1990, 2005, for some history, or Imbens and Rubin, 2006, for relationships to the history of causal inference in economics).

To illustrate, take a look at the specific case in Display 2, where Doctor A is a high prescribing doctor, writing 10 scripts at time 1, and 15 scripts at time 2, whether visited in between or not. Clearly, even though Doctor A writes a large number of scripts, there is no ROI to visit this doctor (for simplicity we are ignoring the cost of a detail, but that is simply a known constant). In contrast, take a look at Display 3, where Doctor B is a low prescribing doctor, writing only one script at time 1, and five or fewer at time 2; yet, Doctor B may be worth visiting, at least much more so than the higher prescribing Doctor A, because a visit to Doctor B will cause an increase in number of scripts from 1 to 5. Whether the four extra scripts, which are caused by the visit, generate a positive ROI for the visit depends on the cost of the visit, the profits from the scripts, and so on.

The point is simply the following: we should make investment decisions based on a comparison of the expected returns when making the investment and when not making the investment: Visit those doctors for whom visiting makes a larger positive difference. Also, allocate company resources to those brands (and those marketing tactics) that provide the greatest marginal positive impact to the company. Great advice (like "buy low, sell high"), but how do we do this in practice?

## 3. THE ESTIMATION OF CAUSAL EFFECTS

The gold standard for the estimation of causal effects is to conduct randomized experiments, such as clinical trials, which are essentially required by the FDA (U.S. Food and Drug Administration) before approving a drug. An alternative, and one which is sometimes acceptable, even to the FDA, is to design and carefully execute an observational study (a nonrandomized design).

Causal effect estimation is not the simple prediction of future events from past events, although these activities can play a role in addressing causal

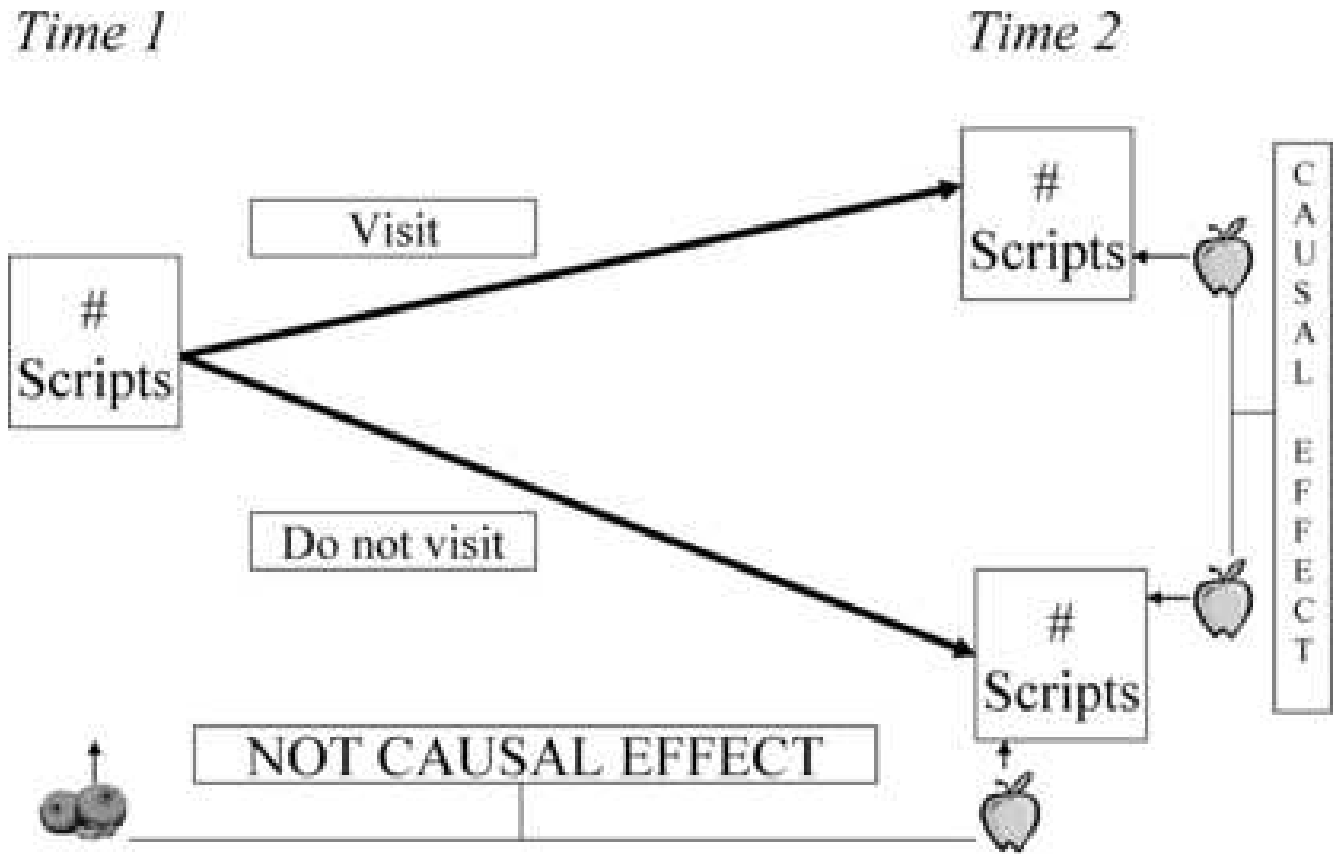

DISPLAY 1. *Causal effect vs. prediction for a particular doctor.*

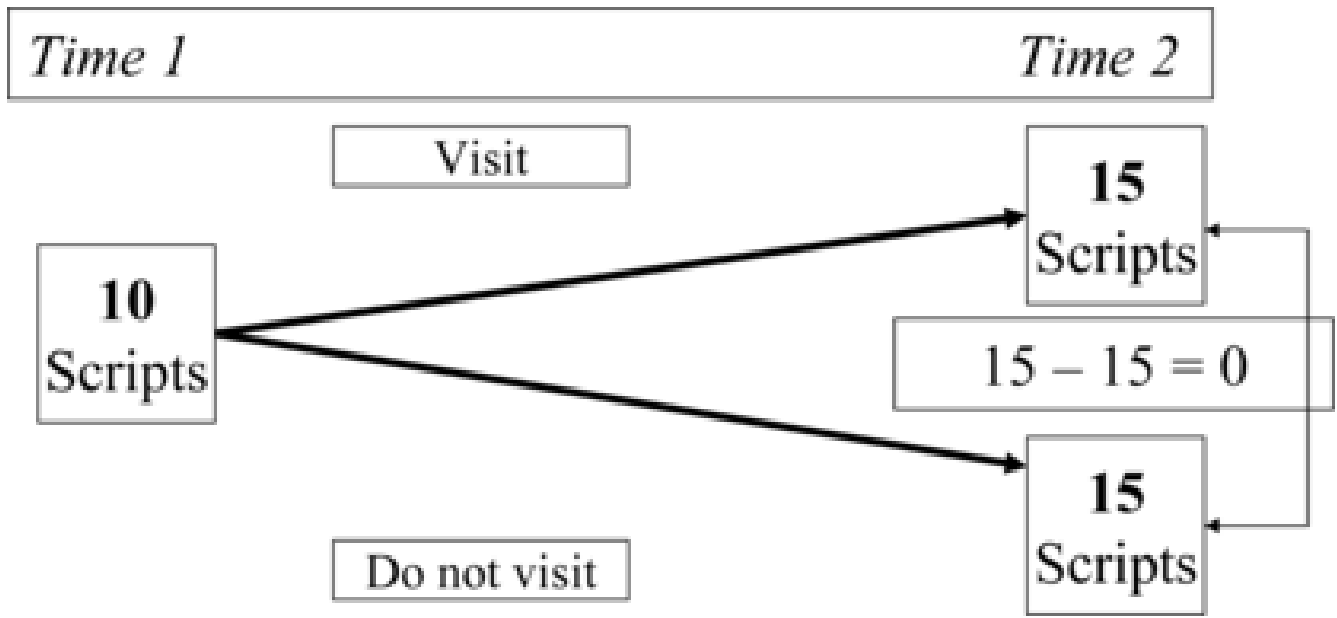

DISPLAY 2. *Example: Doctor A. Temptation is to confuse "prediction" with "causal effect estimation." Example: Visit high prescribing doctor. Waste of money to visit this doctor. Intervention has no causal effect. Here Causal Effect = 0.*

questions. Thus, causal effect estimation is not generally accomplished by: regression, data mining, neural nets, CART, support vector machines, random forests, and so on. Although such techniques can be helpful, none is central, and they can be especially helpful after causal effects of the intervention for each unit (e.g., of the visit for each doctor) have been estimated because it may often be of interest to classify doctors into subgroups based on background variables describing types of doctors, where the subgroups differ by the expected size of their causal effects; this would help future targeting efforts—more on this in Section 9.

So, specifically, how should we think about causal effect estimation from real data? This is easy to describe in principle from the hypothetical database depicted in Display 4. Each row in the matrix displayed there represents one unit (e.g., one doctor), and the columns represent the measurements on them: number of scripts written at time 1, background variables such as age, sex, race, place of doctoral degree, years of practice, type of practice, the number of scripts written by time 2 if visited between time 1 and time 2, the number of scripts written by time 2 if not visited, and the causal effect of being visited—the difference between the latter two. The checked boxes represent observed data values, and the question marks represent unobserved or missing values; the causal effects are all either: (1) a check minus a missing, or (2) a missing minus a check, and both (1) and (2) are effectively missing, and so the entire column of causal effects is always effectively missing.

We describe the process for estimating causal effects as "cloning" for causal effect estimation. Essentially, for each doctor who was visited, we seek a "matching" doctor, a "clone," who was not visited, and we use that doctor's observed outcome (i.e., number of scripts at time 2) to fill in for the first doctor's missing outcome. Similarly, for each doctor who was not visited, we seek a matching doctor, or clone, who was visited, and we use that doctor's observed outcome to fill in for the first doctor's missing outcome. If no matching doctor can be found, the database at hand cannot support causal effect conclusions, at least not without relying on assumptions outside the database (e.g., time-2 scripts without the visit equal time-1 scripts).

Display 5 is a simple reworking of Display 4 where all the missing "?" in Display 4 have been replaced by "!" to indicate that the missing values have been "found" (really "imputed") by using the clones. Of course, finding exact clones for everybody in any real problem is essentially impossible, and yet the conceptual foundations of the above cloning approach rely on using all background variables used in making the decisions to visit one doctor and not visit another. Thus, the pressure to collect many such

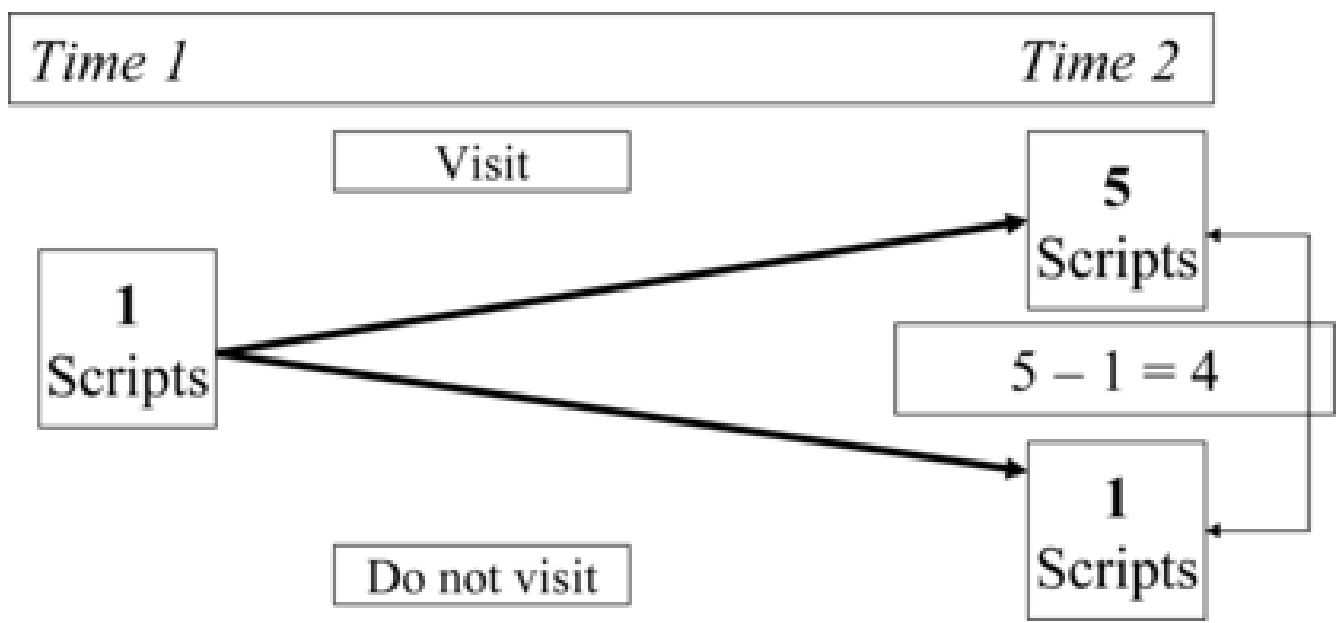

DISPLAY 3. *Example: Doctor B. It is much better to visit this doctor. The investment pays off. Here Causal Effect* $= 4$. *Pays to visit to increase business.*

| Number of scripts at time 1 | Background Characteristics | Number of scripts at time 2 if visited | Number of scripts at time 2 if not visited | Causal effect difference |
|---|---|---|---|---|
| ✓ | ✓ | ✓ | ? | ✓ - ? |
| ✓ | ✓ | ✓ | ? | ✓ - ? |
| ✓ | ✓ | ? | ✓ | ? - ✓ |
| ✓ | ✓ | ? | ✓ | ? - ✓ |

Legend: ✓ Observed, ? Unobserved

Display 4. *Database.*

background variables is great, which, in turn, makes it essentially impossible to find exact clones for anyone.

The key idea for simplification is to use "propensity score" technology (Rosenbaum and Rubin, 1983). This approach allows all the covariates to be reduced to a single covariate. This single covariate, the propensity score, is essentially the probability of being visited as a function of all of the background variables, as estimated from the database. Our actual implementation incorporates other important adjustments and refinements, but the essential idea, and the one being discussed here, is to clone based on the propensity score. For a simple review of some ideas underlying propensity scores, see Rubin (1997, 2006). And for some more recent theoretical work, see Imbens (2000) and Imai and van Dyk (2004); for a couple of hundred thousand other references, just Google "propensity score." Instead of going into details here, we describe our example, and how propensity scores work there.

## 4. ILLUSTRATION: AN ANONYMOUS CASE STUDY

The objective of our actual case study was to produce a target list of doctors based on the estimated effects of a marketing intervention, where the doctors who are thought to provide the best ROI would be visited, at least before the others. That is, we want the doctors ranked by their estimated causal effects (due to a visit) on the number of scripts they would write.

The database consisted of approximately 250,000 doctors in the United States who were active in the medical area of the drug to be promoted. The script data came from an industry standard physician-level prescription database, the sales intervention data came from the company's call reporting system and

DISPLAY 5. *Database after filling in the blanks.*

the doctors' background characteristics came from various other sources; these characteristics included specialty, region of the United States, dates of degrees, and more than a hundred other such variables. The company was (and still is) a top tier U.S. pharmaceutical company.

Display 6 shows the distribution of what is considered to be the most important determinant of whether a rep should visit a doctor: the number of prescriptions of this class of drugs written in the recent past. Notice the rather huge distributional difference between the number of scripts (at time 1) for those who were not visited on the left, and the number of scripts (at time 1) for those who were visited on the right. The doctors who were visited between time 1 and time 2 wrote about 50% more prescriptions per doctor at time 1. Why? Possibly because the salary compensation of the sales reps who visited the doctors was not tied to causing a difference but more to the number of total scripts written by the doctors whom they visited (so they polished the Ferraris!).

But distributional differences between those not visited and those visited are not confined to the number of scripts written at time 1. Display 7 shows the distributions of their specialties (General Practice, Family Practice, Internal Medicine, Endocrinology, OB/Gynecology, Cardiology) for doctors who were not visited (left) and who were visited (right). General Practice was more highly visited, as was Family Practice and Internal Medicine and Endocrinology, whereas OB/GYN was not visited as much. Similar displays could be created for dozens of other background variables.

Display 8 presents a single picture that summarizes, at least to a great extent, all such displays. It presents the distributional comparison for the estimated propensity scores for those not visited (top)

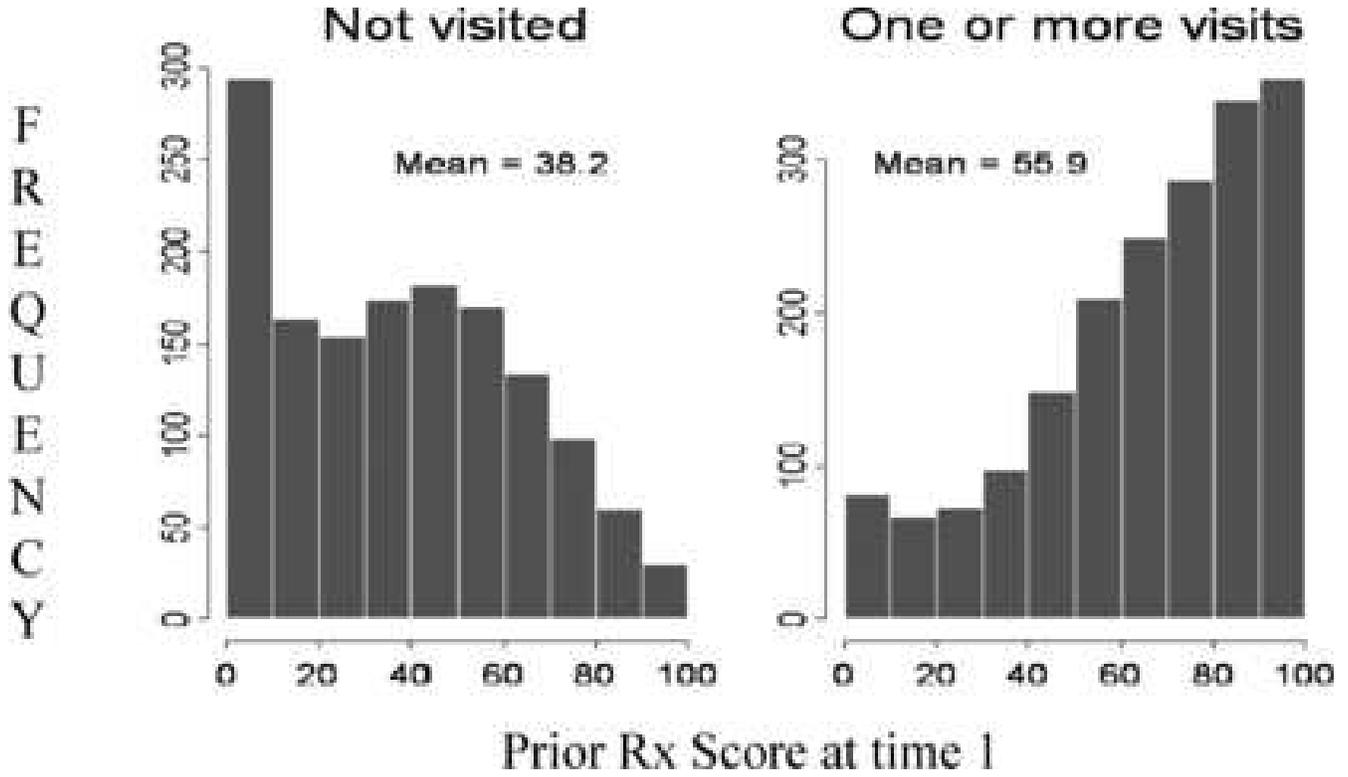

Display 6. *Histograms for background variable:* Prior Rx Score at time 1.

and those visited (bottom). The label "linear propensity score" means that instead of plotting the estimated probability of being visited, say $e$, we plotted $\log[e/(1-e)]$; the reason for doing so is technical, and is discussed briefly in an applied context in Rubin (2001). The propensity score, by construction, is supposed to lead to the worst (i.e., most discrepant) comparison of distributions of any combi-

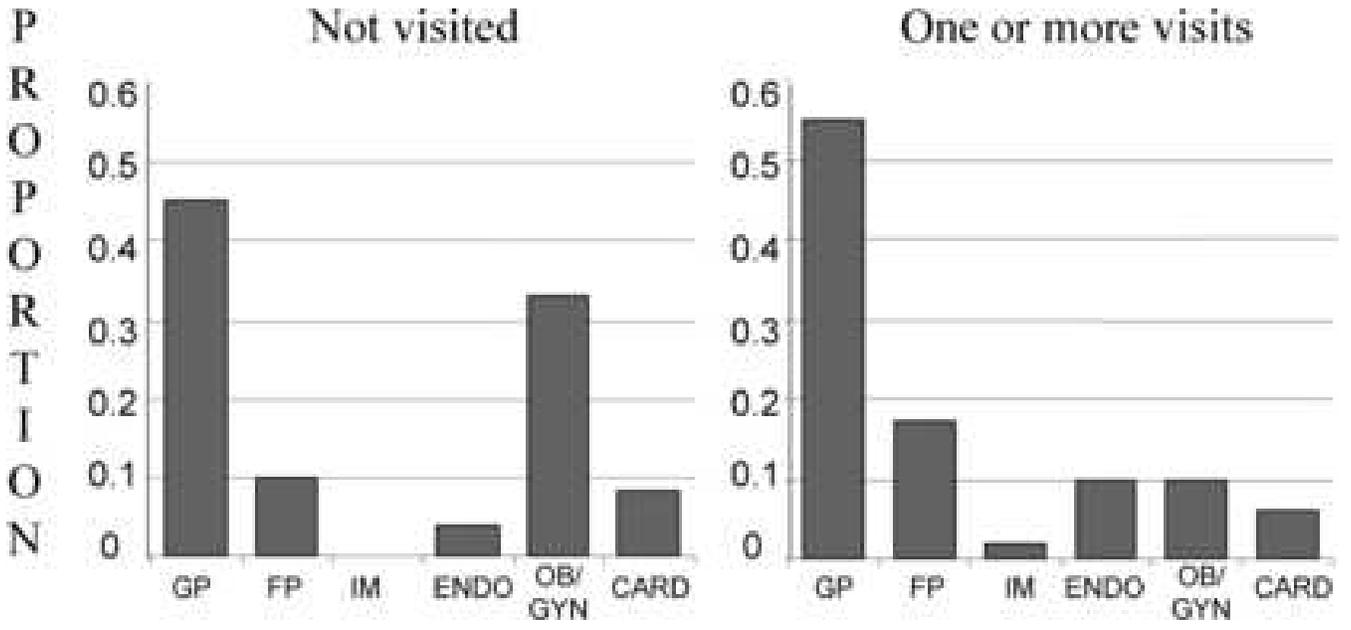

Display 7. *Histograms for background variable:* Specialty.

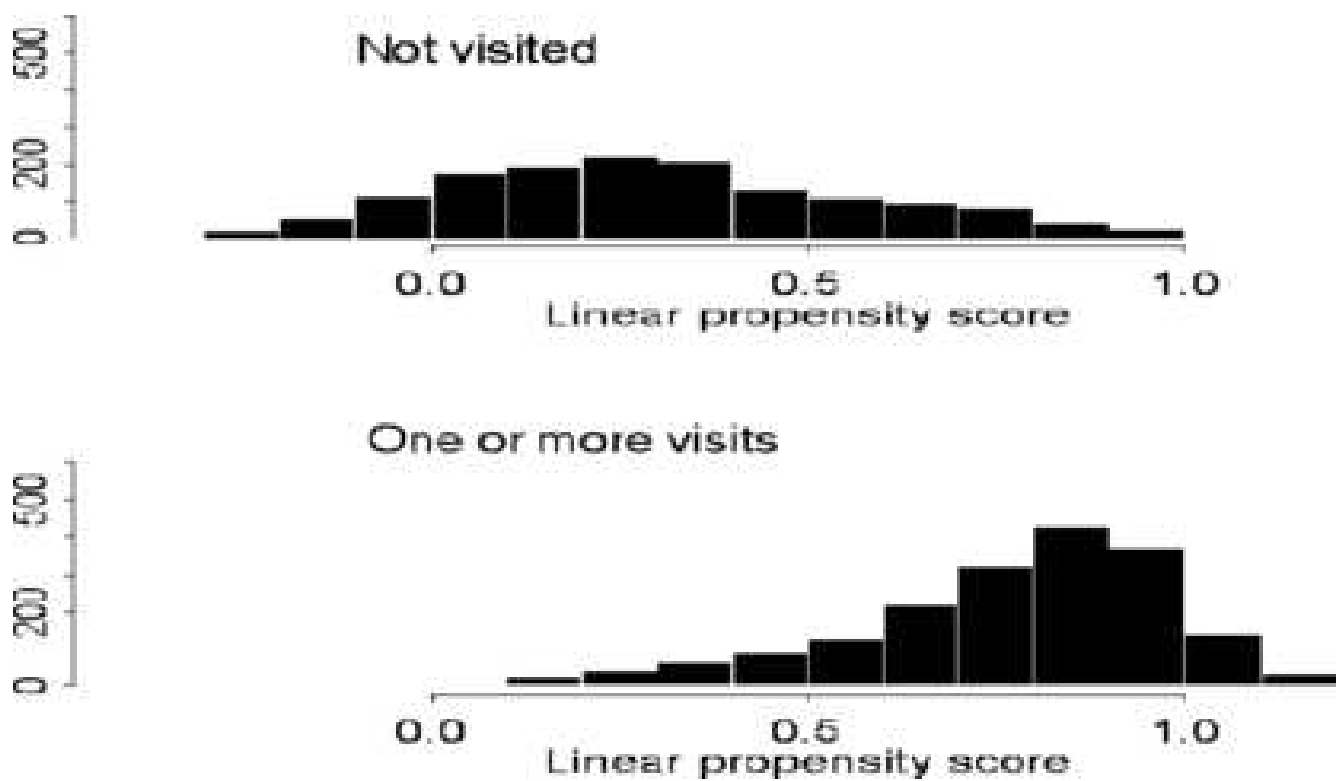

Display 8. *Histograms for summarized background variables:* Propensity Score.

nation of the underlying background variables (here "background variables" includes the scripts written at time 1). In this example, as in many others, the estimated propensity scores are found by running a logistic regression, with visited/not visited as the dependent variable, and all the background variables (including time-1 scripts) as predictors. Any good method for predicting the probability of being visited or not could be used; this is another possible use for statistical software designed for prediction or classification.

The differences revealed in Display 8 between the not-visited and the visited groups in the distributions of their estimated propensity scores are dramatic. Even if we were to allow cloning to be done only using this variable, there are many not-visited doctors without clones (i.e., those with linear propensity scores below about 0.1), and there are a collection of visited doctors without clones (i.e., those with linear propensity scores above 1.0). However, for the other values of the linear propensity score, it appears that there will be some clones, at least if we only consider the propensity score and ignore the other 100+ background variables.

But can it be valid to focus on this one "compound" variable and ignore the others? In a certain sense, yes, it is valid. A theorem due to Rosenbaum and Rubin (1983) implies the following claim in large samples: Suppose we look in one narrow bin of propensity scores, say from 0.5 to 0.6, and we compare the distributions of any background variable for the not-visited and the visited. We will find that the two distributions are essentially the same, even though, for some variables (like the number of visits at time 1 given in Display 6), the overall, unbinned, distributions for the visited and not visited are very different. Perhaps this sounds surprising, but the two distributions for the number of scripts at time 1 in this propensity score bin are given in Display 9. These two distributions in Display 9 are not identical, but they sure are close, much closer than in Display 6. Display 10 provides the analogous distributions in the same bin of propensity scores but for doctor's specialty. Again, they are not perfectly matched, but they are much closer than the overall, unbinned, distributions given in Display 7.

The same thing will be found in this bin for each of the 100+ background variables that were used to estimate the propensity score, as well as for each of the propensity score bins with a reasonable number of doctors who were not visited and who were

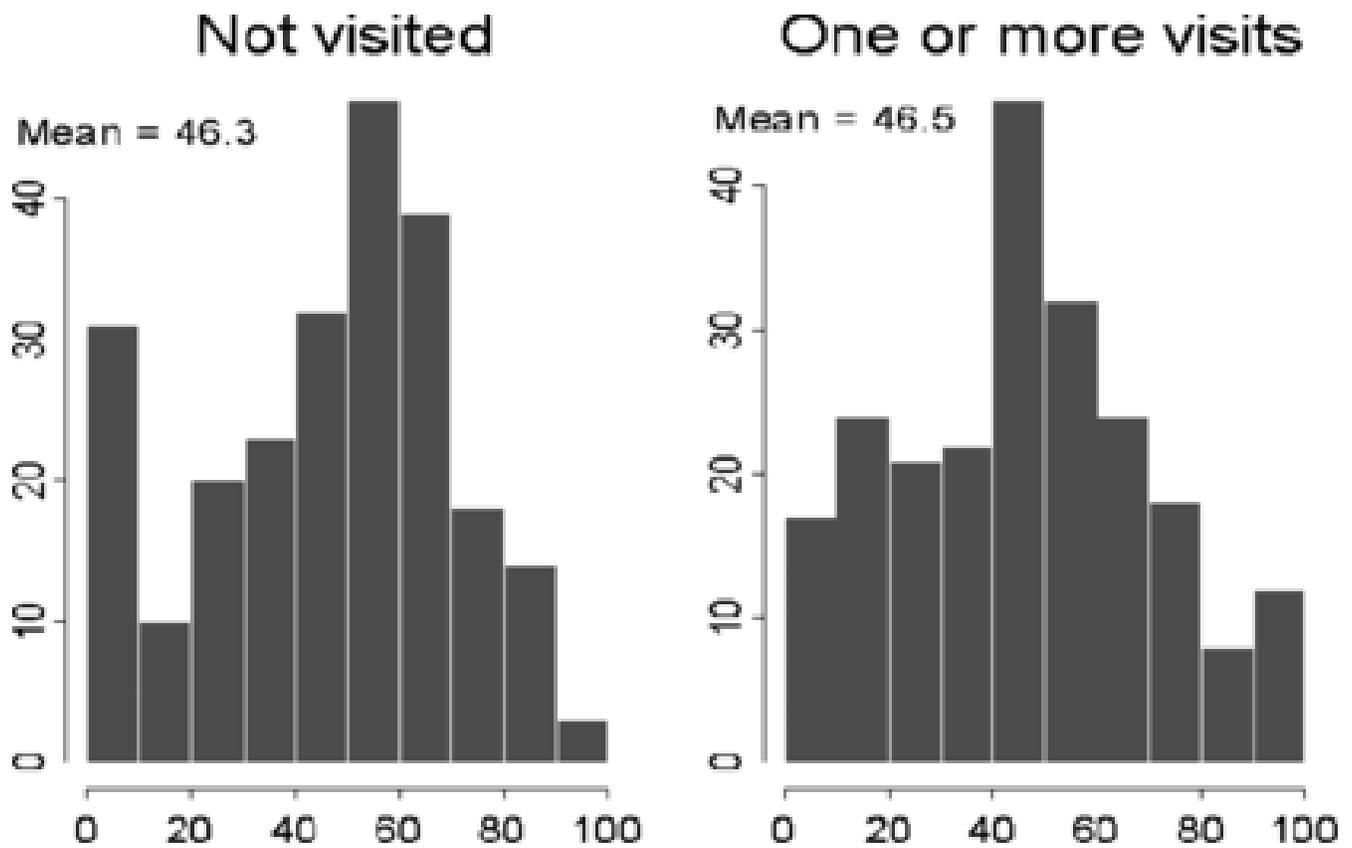

Display 9. *Histograms for a variable in a bin of propensity scores:* Prior Rx Score.

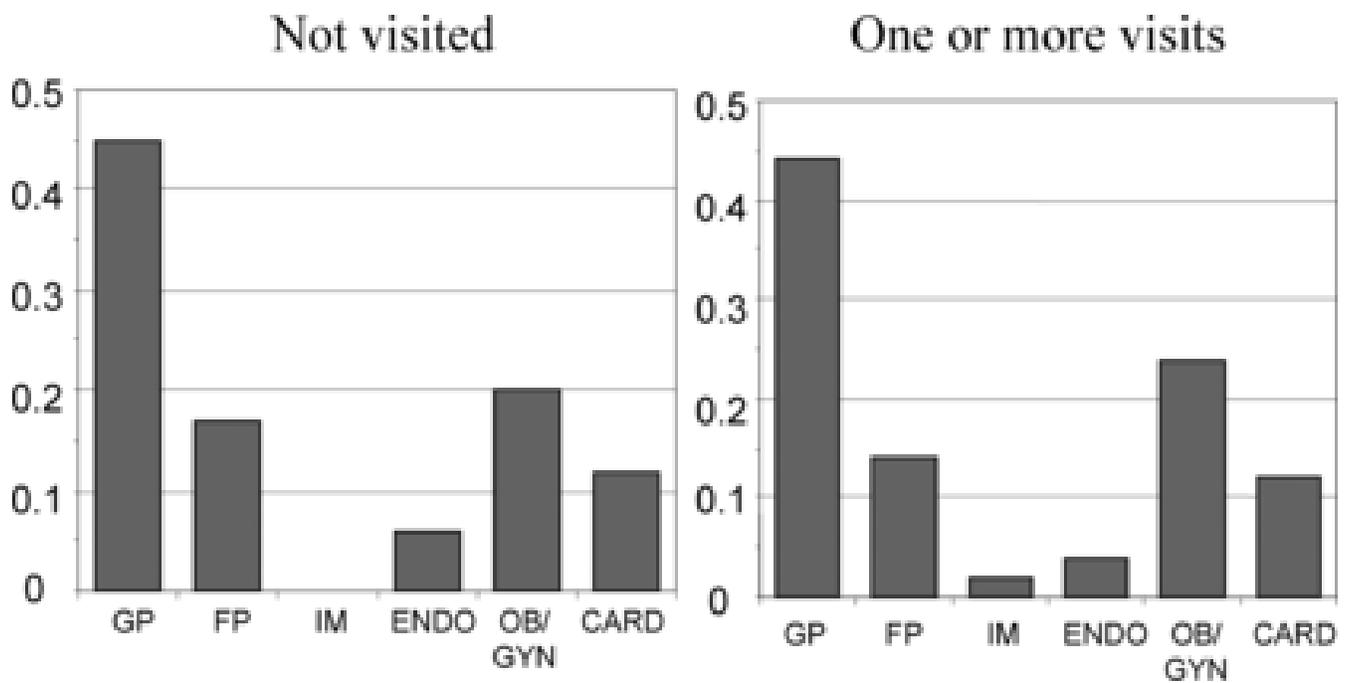

Display 10. *Histograms for a variable in a bin of propensity scores:* Specialty.

visited on which to base the distributional comparisons. The intuition behind this result is the following: let us suppose we have a group of visited and not-visited doctors, all of whom have the identical propensity score (i.e., a very narrow bin, between 0.29999 and 0.30001, essentially all at 0.3 estimated probability of being visited). Then each doctor in this bin has probability 0.3 of being assigned a visit, and so the not-visited and the visited are only randomly different from each other. And randomization assures the same distribution of all background variables, at least in large samples. But because the propensity scores were estimated using only some background variables, this conclusion only holds for the variables used to create it, not for all background variables, as with true randomization.

Of great importance, notice that the cloning takes place using only background information and the information about which of the doctors were visited. No outcome data, that is, no time-2 script data, were available. This is imperative for the honesty of the cloning and is in stark contrast to model-based predictive approaches (e.g., see the discussion in Rubin, 2001, 2005; Imbens and Rubin, 2006); this point is emphasized and discussed further in Section 9.

## 5. DOCTOR-LEVEL CAUSAL EFFECT ESTIMATION IN THE CASE STUDY AND RESULTANT RECOMMENDATIONS

Now our database also had, in addition to whether each was visited or not between times 1 and 2, and the background information, the outcome of interest: the number of scripts at time 2 for each of the doctors. For cloning, this variable was "held in escrow"—not available. After the clones were created, we used those data as follows. Having clones identified for doctors using their propensity scores, we filled in the missing values in Display 4 using the clones' observed values, and thereby created estimated causal effects for each of the doctors who has a clone. As stated earlier, some of the doctors, those with extremely large or small propensity scores, did not have clones, and so did not have estimated causal effects; such estimates could have been created, but this would rely on more assumptions, and that topic is beyond the scope of this very simple introduction to the ideas, and is only briefly addressed in Section 9.

Having created estimated causal effects for the visits, which were scheduled between time 1 and time 2, on scripts written at time 2, we then focused on the question of deciding which of the previously unvisited doctors to visit between now (after time 2) and time 3, in the future, to maximize ROI. The answer was: order the not-yet-visited doctors by their estimated causal effects of a visit between time 1 and time 2 on scripts at time 2. We have no data from time 3, the future, so all that we could do is hope that the causal effects in the future will be like those in the past. That is, for the not-yet-visited doctors, we estimated the causal effects of a visit between times 2 and 3 by the estimated causal effect of a visit between time 1 and time 2, which was simply their estimated number of time-2 scripts if visited (based on the cloning) minus their observed number of scripts at time 2 when not visited. (As mentioned earlier, there were some adjustments and refinements that we used in the actual application briefly discussed in Section 9, but the basic idea is not changed.) And then the not-yet-visited doctors were ordered by their estimated causal effects.

The obvious recommendation was to visit those doctors at the top of the list of estimated causal effects: the top 10% first, then the next 10%, and so on. The company's then-current practice was to recommend visiting doctors according to a "standard decile" list. This list was based on regression statistics, which predicted the number of scripts written at time 2 from scripts written at time 1 and background variables. And the company's recommendation was to visit the top 10% first, then the next 10%, and so on. The essential difference between these two recommendations was how the list was created: by "cloning for causal effect estimation" versus regression prediction. Only in very special circumstances would the two lists generate the same or even similar orderings.

We were also asked by the company to estimate the result of using their list, which was more or less being implemented, with a lot of local decision-making that did not follow it, and what would have happened if instead they had used our list to order the doctors.

## 6. VALIDATION OF ORDERING BY ESTIMATED CAUSAL EFFECTS

Among the previously (at time 2) not-visited doctors, some were visited between time 2 and time 3, and some were not visited. As stated above, the sales reps were not very adherent to the lists they were supposed to follow, possibly because of their own assessments of who was going to be a big script writer in the near future, possibly because of personal connections, and so on. The company had the data on the number of scripts written in a November and the following January, and they wanted to estimate the relative effectiveness of the two lists for predicting the effect of a December detail, where they considered the change in the number of scripts written from November to January to be the estimated causal effect of the visit for each of these previously unvisited doctors.

In particular, we were to order all visited-in-December doctors according to both lists, then compare the average estimated causal effect in the top decile of our list and the average estimated causal effect in the top decile of their traditional list. And we were to continue down both lists, comparing decile by decile. The results are shown in Display 11.

Although all estimated average causal effects were small (to be fair, the time involved to see any causal effects was very short), our top decile was much superior to their traditional top decile. In fact, their list revealed a negative estimated causal effect of a December visit on scripts for the doctors in their top decile, whereas the doctors in the top decile of our list revealed a positive estimated causal effect. Our corresponding seventh, eighth and ninth deciles also included doctors with positive and superior causal effects; deciles 2 through 6 were noisy, but suggested positive causal effects for doctors in our deciles. Of substantial interest, the doctors for whom we predicted the worst effect of a December visit, our decile 1 group, had the worst causal effects—quite negative. In contrast, the doctors in the bottom decile for the traditional list turned out to have had a larger positive estimated causal effect of the December visit than any other decile from the traditional list—the exact opposite of what the company hoped! Even more dramatic results were obtained for the causal effects of 8+ details over more months—see Display 12.

The conclusion from this validation phase seems fairly clear: Use the estimated causal effects to create targeting lists; do not use the traditional list based on predicted scripts. This policy change should result in a better return on investment. Of course, there are various institutional barriers to any policy change that is dramatic, such as following this recommendation. For example, the way sales reps are compensated may have to be altered so that they do not get rewarded for "polishing Ferraris" to make them faster, to use the analogy from earlier.

Targeting individual doctors is just one of the many tasks that can be performed with such causal effect estimation. For example, once individual doctor-level causal effects are estimated, all of the usual prediction and classification methods can be highly appropriate for the prediction of types or subgroups of doctors outside the original database who are worth targeting because similar types in the database are predicted to have large estimated causal effects. This will be discussed in Section 9 in more depth, after brief discussions of some details of the methodology and some potential applications in e-commerce.

## 7. MORE DETAILS ABOUT THE METHOD

This section provides a more explicit roadmap for the practitioner regarding the steps involved in such an analysis. This roadmap, Display 13, is presented by way of a flow diagram with accompanying brief discussion of each of seven steps that can be used to define the approach.

1. Assemble data sources. Though self-evident in its necessity, this step may well be one of the most time-consuming practically, because a data source that has a rich collection of covariates is frequently one that has been assembled from a variety of sources. In the example here, it would comprise prescribing data, physician demographic data, elements from the *call database* (e.g., the data source that contains information on the rep's details), geographic information (rural/urban, e.g.) and HMO coverage indicators.

Of course the individual units in such data sources need not be individual doctors or even people. These units could be households; they could be medical practices; or they could even be cities, with radio advertisements serving as the intervention. Another possible unit is that of a retailing location (Singh, Hansen and Blattberg, 2006).

2. Establish the definition of the intervention or "promotional event." For example, we could formulate the intervention as zero details versus one detail, or we could formulate it as zero details versus one or more details. It is up to the practitioner to decide which question is the most relevant. A common notation for the intervention is to denote it by the indicator $Z$, where $Z = 1$ for the units that received the intervention and $Z = 0$ for those that did not.

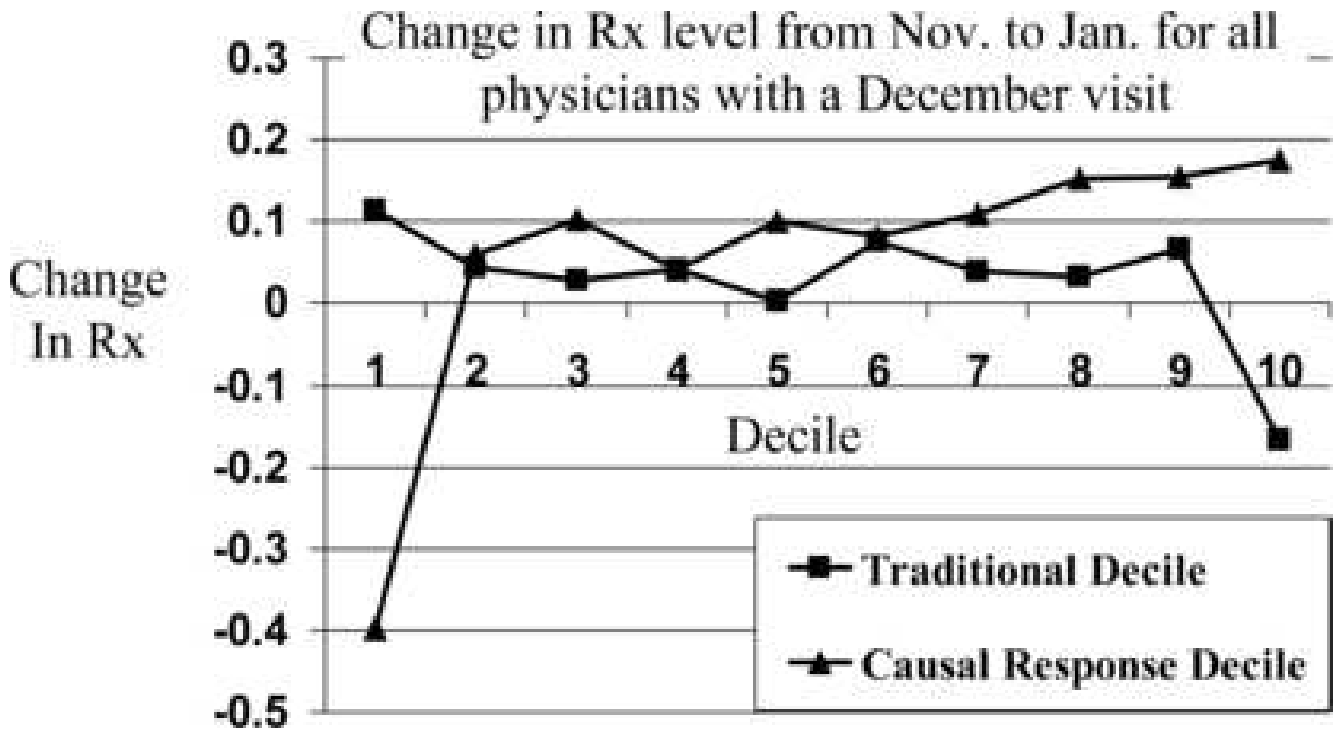

Display 11. *The estimated effect of a single visit.*

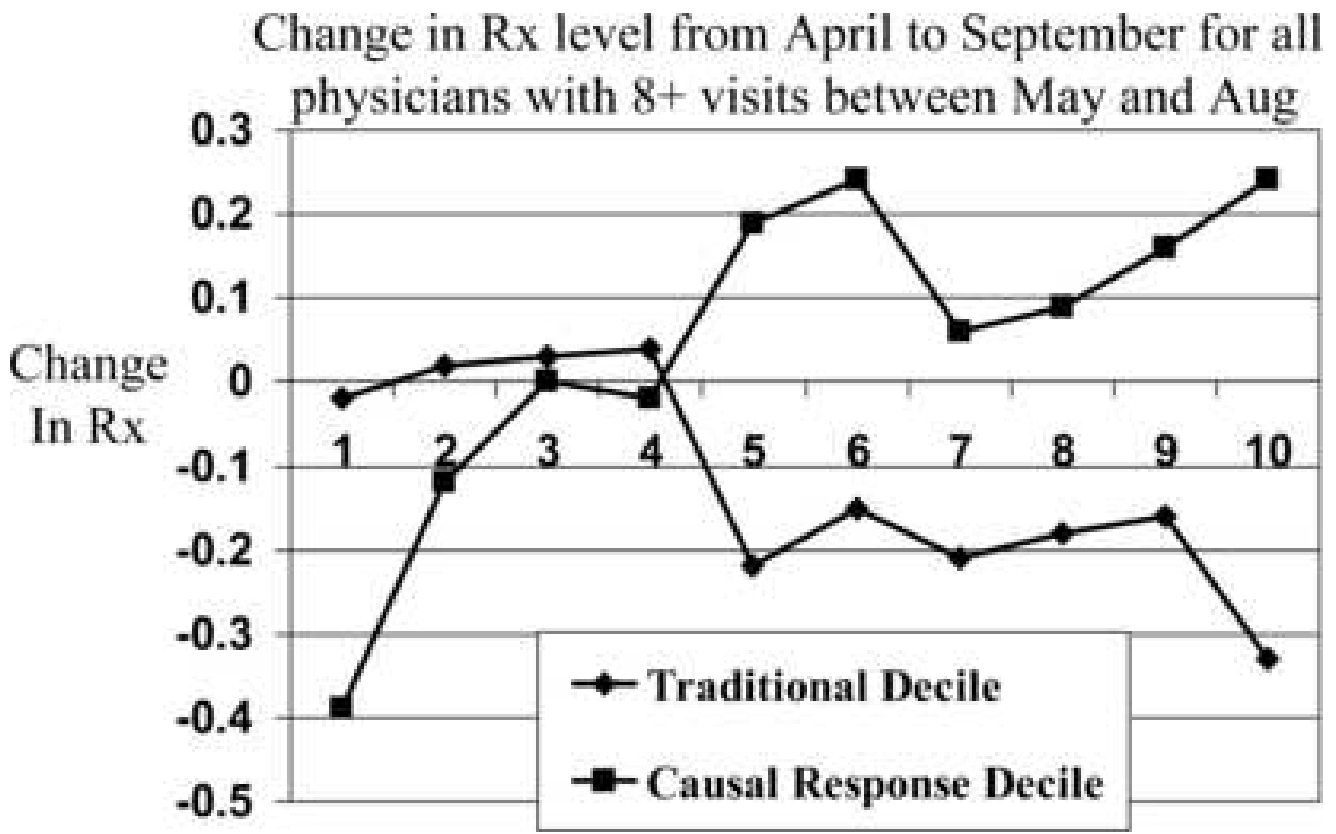

Display 12. *The estimated effect of multiple visits.*

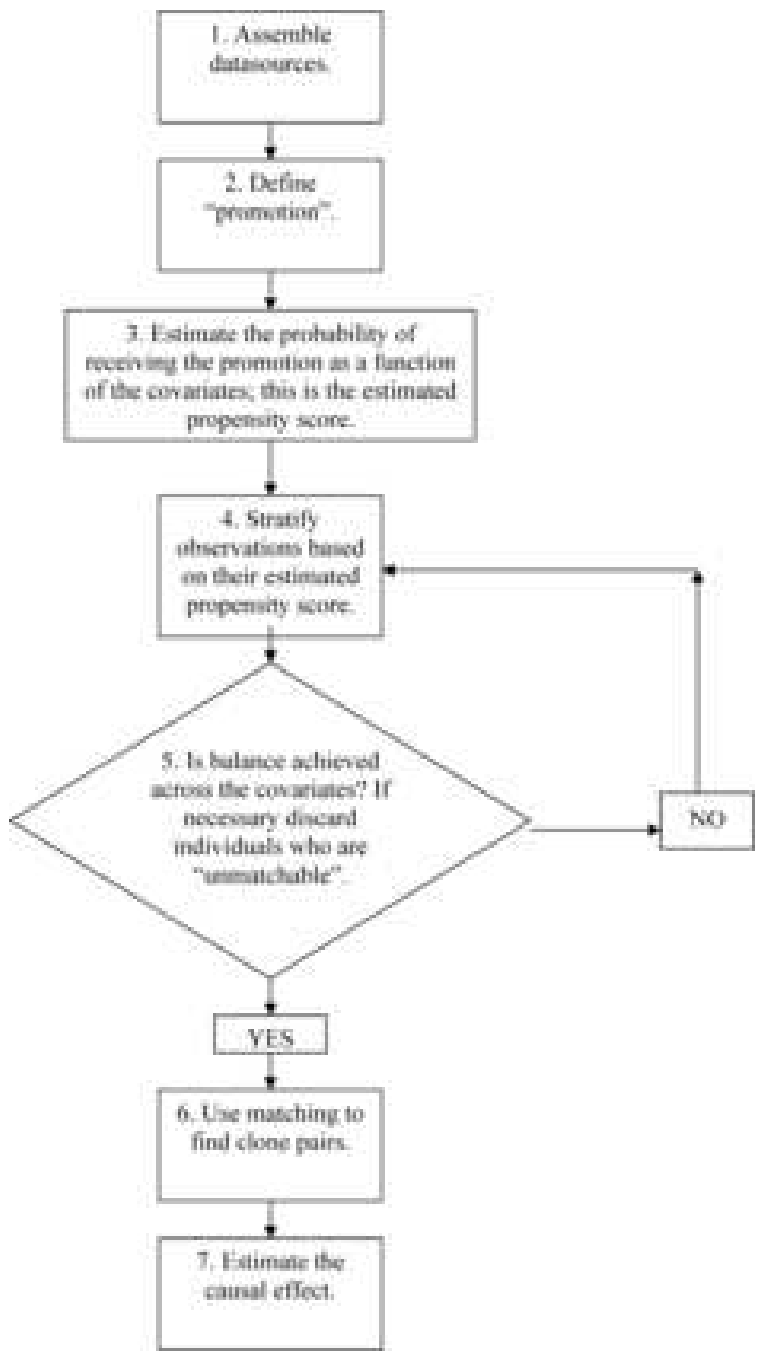

DISPLAY 13. *A roadmap for implementing the basic propensity score methodology.*

3. Estimate the unit-level probabilities of receiving the intervention as a function of the unit's values of covariates; that is, estimate $e = P(Z = 1|X)$ where $X$ denotes the covariates. These probabilities may be estimated by a variety of methods, the most common probably being logistic regression, but any predictive modeling methodology could be appropriate here.

4. Stratify (or bin) the units into groups based on their estimated values of $e$, that is, their estimated propensity score.

5. Run diagnostic tests (see Rosenbaum and Rubin, 1984; Rubin, 2001, as illustrated in Displays 9 and 10) to ensure balance has been met. Balance is used in the sense that the distribution of the covariates within the propensity score bins should be the same regardless of whether $Z = 0$ or $Z = 1$. One of the common misconceptions with the methodology is the belief that the individual clone pairs have to be very similar on all characteristics, as in the common understanding of a matched pair, but that is in fact not the case. They need only be similar in this distributional sense as just described when global causal effects are to be estimated. Similarly, if more refined questions are asked within specific subgroups, then distributional balance in the covariates needs to be achieved within the subgroup.

If balance has not been achieved, refine the model for the propensity scores until it is satisfactory, that is, return to step 3. It is often the case that some of the estimated propensity scores are so extreme in either the positive or negative direction that no close match can be found. In this situation, we remove these units from the analysis with the acknowledgment that the data set is uninformative with regard to estimating their causal effects without relying on extraneous assumptions. In an ideal world, we would recommend this group for experimentation, though practically the units with very high estimated propensity scores tend to be high value with respect to revenue, so that there will likely be a strong reluctance to "experiment" with this group.

6. Use the matching to find clone pairs. Though straightforward in concept, there are many available algorithms for this step, including "greedy" and "optimal" methods (e.g., see Rosenbaum, 1989, and Gu and Rosenbaum, 1993). In addition, there may be benefits achieved in moving beyond one-to-one matching, and in the work presented here we used one to ten matching to reduce the variability of the casual effect estimates, possibly at the expense of a little bias (e.g., for some discussion, see Ming and Rosenbaum, 2000, and Rubin and Thomas, 2000).

7. Estimate the casual effect. The elegance of the approach is really revealed in the final estimation step where the outcome variable is simply differenced between the elements of a clone pair. This provides an estimate of an individual level causal effect. These individual level effects may be aggregated within any subgroup of interest, and in the pharmaceutical example presented here, one could, for instance, be particularly interested in the efficacy of regional sales forces, or the deciles of a target list. We take one additional opportunity to point out that this is the only place in the entire analysis where the outcome variable is used (the Y in "escrow" idea again).

## 8. OPPORTUNITIES FOR APPLICATIONS IN E-COMMERCE

Propensity score methods can be applicable whenever there is a desire to estimate the impact of an intervention, and it is particularly relevant when the intervention is not applied on a randomized basis but we think we have the major background variables that influence which treatment is received.

To illustrate with an e-commerce example, take an online store with an e-mail membership database. Some of these members have opted-in to receive a monthly e-mail newsletter, but the company does not know whether the resources they put into producing the newsletter make financial sense. It is clear that the newsletter has not been randomized to members, and indeed, those individuals who opted to receive the newsletter may be systematically different from those who did not.

The propensity score approach to this problem would be to use a set of covariates describing all the members, which should almost certainly include prior purchasing behavior, and create a model for the probability they received the newsletter. With this estimation completed, that is, with the propensity scores estimated, the matching can proceed and the causal effect of the newsletter estimated by the average difference in revenue between the clones, those receiving the newsletter and those not receiving the newsletter.

Another similar example would be to find the ROI of a free shipping initiative, where again individuals opt-in to receive the free shipping but pay a modest yearly fee for it. Again, the "treatment"—free shipping—is not randomized to respondents, so a naive comparison of the purchasing behavior of those who receive it to those who do not would not lead to a valid ROI estimate.

A subsequent use ideally suited to e-commerce could follow from this previous activity. If a predictive model for the individual-level causal effects of the free shipping promotion themselves were created, then this model could be used to determine who should be targeted with the free shipping; essentially we would not want to offer the free shipping to someone who would have spent the money anyway, but rather we would like to target the free shipping toward those for whom it is most effective, where effective is understood in the sense of those with the largest causal effects of the free shipping.

A particularly interesting opportunity would occur if activity on the website itself could be used to generate some of the background covariates in this predictive model for causal effects. To illustrate with an example that is based on our experience, consider a pharmaceutical company with a prescription drug and a consumer-orientated website providing information on the drug. The offering of a coupon to encourage trial use of the drug is common on such sites. But to whom should the coupon be offered and for how much? Ideally we would offer it to the more "valuable" customers. With a predictive model for the causal effect of the coupon using covariates based on website behavior as measured through the click-stream log (see Bucklin et al., 2002; Montgomery et al., 2004), the company could have an on-the-fly method for generating such offers. Indeed, we have found, in previous proprietary work, that page trails and specific pages visited on a website were related to potential customer value, albeit self-reported through intent to purchase. However, the potential for the methodology to be used in the targeting arena is certainly strong with a predictive model for causal effects available.

When individual causal effects are a focus, refinements to the estimates based solely on clones can be highly useful. We now turn to a brief discussion of such refinements which use statistical models to help create smoother unit-level causal effect estimates.

## 9. THE ROLE OF MODELS WHEN ESTIMATING CAUSAL EFFECTS

Of substantial importance, the propensity score approach to causal inference in business advocated here has focused on the theme that the design of an observational study should parallel the design of a randomized experiment. That is, our propensity score approach is accomplished without any access to outcome data (the Y-variable in escrow paradigm), and it seeks to create samples of units who were exposed to the intervention and who were not exposed to the intervention that are as similar as possible to each other pre-intervention. This approach is in stark contrast to the standard approach in economics and statistics based on the fitting and refitting of straight lines or similar models, where the estimated causal effect is given by the estimated coefficient of an indicator variable for exposure to the intervention. This approach has been known for many years to fail in general (e.g., Cochran and Rubin, 1973; LaLonde, 1986). In these same situations, matching alone or in combination with modeling has been successful (e.g., see Dehejia and Wahba, 1999; Rubin and Thomas, 2000, and their references).

A real problem with the traditional approach to causal inference is that the estimated answers are constantly being seen and altered as models are fit and refit. Also, because typically the intervention and control groups are far apart, for example, as revealed by Displays 6, 7 and 8 in our example, different models will effectively imply different extrapolations, with possibly wildly differing answers from which the analyst can pick and choose. This picking and choosing is often done by someone who is effectively in a "conflict of interest" position because there is often knowledge of what the client would like to see concerning the intervention being evaluated. This is no way to achieve honest answers, even if the data analysts are honest—the unconscious pressures to find one sort of "appealing" answer are often tremendous.

The advantage to the modeling approach is that when the models are appropriate, the resulting causal estimates can be more efficient and precise than the propensity score cloned estimates. This is no reason to give up the objectivity of the propensity score approach, but it does suggest using the modeling approach to refine the causal estimates based on the matched samples. In fact, we have applied such refinements with success to obtain more reliable unit-level causal effects in our own work.

Additionally, once we have reliably estimated individual causal effects, with or without the aid of models it can be highly useful to try to predict these estimated causal effects from background variables that are available in the samples used to estimate the causal effects, as well as in future groups of units that we may wish to target for the intervention. Obviously, the type of unit that is predicted to have a large estimated causal effect is the one to target first. Also, such predictions can be very important for allocating resources. For such allocations, standard statistical prediction models can be entirely appropriate, and their use does not interfere with the objective and honest estimation of the basic unit-level causal effects based, fundamentally, on the creation of the clones.

We hope that this simple introduction has stimulated more appreciation for the relevance and importance of taking the estimation of causal effects seriously, and how it generally differs from simple prediction.

## ACKNOWLEDGMENT

We are indebted to our colleague Scott Nass for presenting the substantive problem to us, and discussions regarding the finer points of pharmaceutical detailing.